\documentclass{article}
\usepackage{hyperref}
\usepackage{amssymb}
\usepackage{amsmath}
\usepackage{amsthm}
\usepackage[nice]{nicefrac}
\usepackage{xfrac}
\usepackage{mathtools}
\usepackage[numbers]{natbib}
\usepackage{optidef}
\usepackage{authblk}
\usepackage{graphicx}
\usepackage[misc]{ifsym}
\usepackage{bbding}
\usepackage{lineno}
\newcommand{\upr}{\overline{P}}
\newcommand{\lpr}{\underline{P}}

\newcommand{\lpe}{\underline{E}}

\newcommand{\rset}{\mathbb{R}}
\newcommand{\natset}{\mathbb{N}}

\newcommand{\zset}{\mathbb{Z}}

\newcommand{\lset}{\mathcal{L}}

\newcommand{\naturals}{\mbox{I\hspace{-.1cm}N}}
\newcommand{\prt}{I\dsn P}

\newcommand{\nset}{\naturals}

\newcommand{\dset}{\mathcal{D}}
\newcommand{\aset}{\mathcal{A}}

\newcommand{\mset}{\mathcal{M}}

\newcommand{\asetpa}{\aset(\prt)}

\newcommand{\dsn}{\!\!}

\newcommand{\comment}[1]{}
\newcommand{\nega}[1]{{#1}^\mathsf{c}}

\newcommand{\lvx}{\underline{V}(X)}
\newcommand{\uvx}{\overline{V}(X)}

{\left\lbrace\begin{array}{@{}l@{}}}%
	{\end{array}\right.}

\newtheorem{theorem}{Theorem}[section]
\newtheorem{proposition}{Proposition}[section]
\newtheorem{lemma}{Lemma}[section]
\newtheorem{corollary}{Corollary}[section]
\newtheorem{remark}{Remark}[section]
\newtheorem{definition}{Definition}[section]
\newtheorem{example}{Example}[section]

\DeclareMathSymbol{\mh}{\mathord}{operators}{`\-}

\begin{document}

\title{Jensen's and Cantelli's Inequalities with Imprecise Previsions}
%
%

%
%


%

\author[1]{Renato Pelessoni\thanks{renato.pelessoni@deams.units.it}}
\author[1]{Paolo Vicig\thanks{paolo.vicig@deams.units.it}}
\affil[1]{DEAMS ``B. de Finetti''\\
	University of Trieste\\
	Piazzale Europa~1\\
	I-34127 Trieste\\
	Italy}

\renewcommand\Authands{ and }

\maketitle

%
\begin{abstract}
We investigate how basic probability inequalities can be extended to an imprecise framework, where (precise) probabilities and expectations are replaced by imprecise probabilities and lower/upper previsions.
We focus on inequalities giving information on a single bounded random variable $X$, considering either convex/concave functions of $X$ (Jensen's inequalities) or one-sided bounds such as $(X\geq c)$ or $(X\leq c)$ (Markov's and Cantelli's inequalities).
As for the consistency of the relevant imprecise uncertainty measures, our analysis considers coherence as well as weaker requirements, notably $2$-coherence, which proves to be often sufficient.
Jensen-like inequalities are introduced, as well as a generalisation of a recent improvement to Jensen's inequality.
Some of their applications are proposed: extensions of Lyapunov's inequality and inferential problems.
After discussing upper and lower Markov's inequalities, Cantelli-like inequalities are proven with different degrees of consistency for the related lower/upper previsions.
In the case of coherent imprecise previsions, the corresponding Cantelli's inequalities make use of Walley's lower and upper variances, generally ensuring better bounds.

\end{abstract}

\section*{Acknowledgement}
*NOTICE: This is the authors' version of a work that was accepted for publication in Fuzzy Sets and Systems. Changes resulting from the publishing process, such as peer review, editing, corrections, structural formatting, and other quality control mechanisms may not be reflected in this document. Changes may have been made to this work since it was submitted for publication. A definitive version was subsequently published in Fuzzy Sets and Systems, 

https://dx.doi.org/10.1016/j.fss.2022.06.021 $\copyright$ Copyright Elsevier

\vspace{0.3cm}
$\copyright$ 2022. This manuscript version is made available under the CC-BY-NC-ND 4.0 license http://creativecommons.org/licenses/by-nc-nd/4.0/

\begin{center}
	\includegraphics[width=2cm]{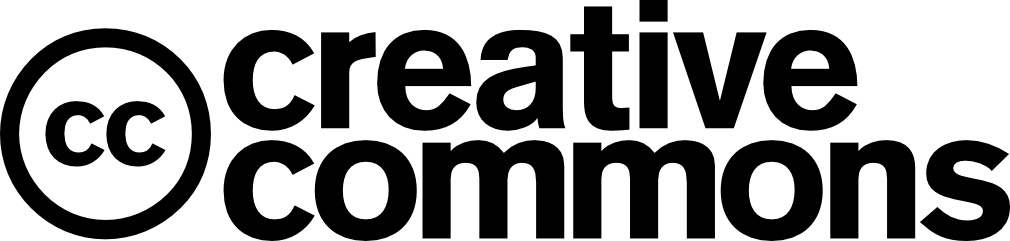}
	\includegraphics[width=2cm]{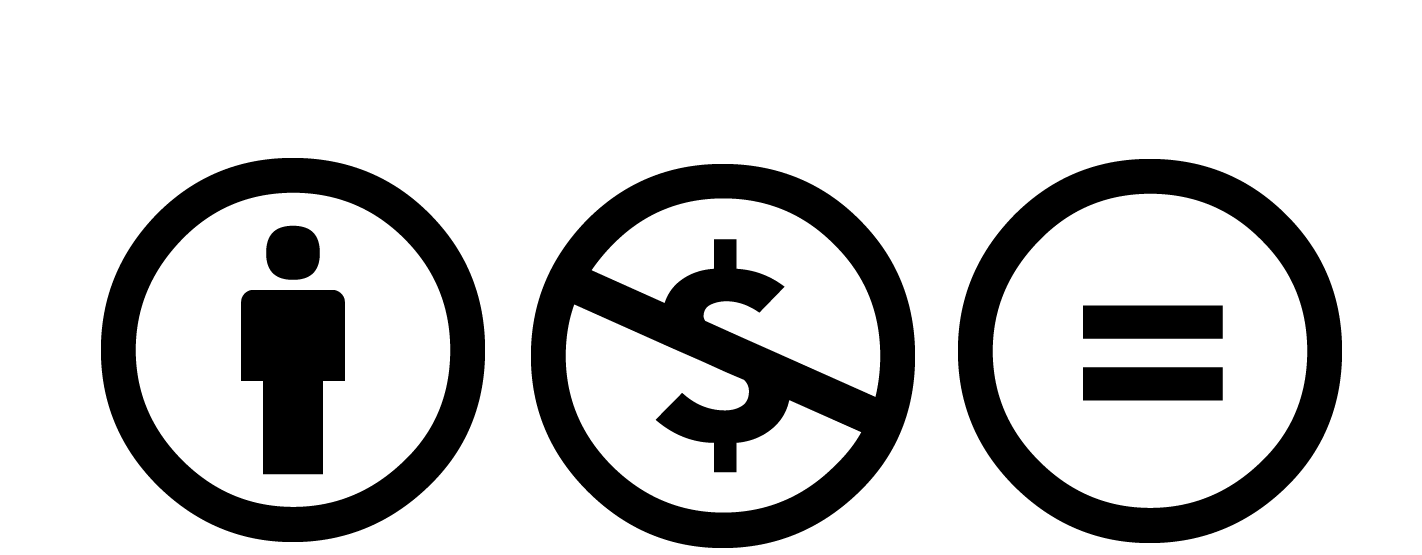}
\end{center}


%
%
%
\section{Introduction}
\label{sec:introduction}
Probability theory is rich of a number of more or less famous inequalities, some of which, like Markov's inequality, usually presented already in introductory courses.
Inequalities may serve several purposes, and a common one is offering bounds to uncertainty evaluations of a random number $X$.
The bounds are often independent of the effective distribution of $X$, hence are especially useful when its knowledge is vague or imprecise.
Thus, a natural question is: how are well known probability inequalities modified in the more general frame of Imprecise Probability theory?
The theory deals explicitly with imprecise evaluations: for instance, lower/upper previsions replace and generalise expectations \cite{walley_statistical_1991}.
A specific feature of lower/upper previsions is that they may be defined on \emph{arbitrary} sets of \emph{gambles} (i.e. bounded random variables). This matches well with practical situations, where effective uncertainty knowledge may be limited to just a few gambles.
Another relevant feature is that more consistency criteria for lower/upper previsions, of different strengths, have been proposed.
\emph{Coherence} \cite{walley_statistical_1991} is the most investigated and closest to the precise probability case.
Weaker notions are \emph{$2$-coherence}
\cite{Pelessoni_Vicig_2016, Pelessoni_Vicig_2017,walley_statistical_1991},
the condition of avoiding sure loss \cite[Section~2.4]{walley_statistical_1991} and others.
Hence, a second interesting question arises: if a probability inequality can be extended in some way to imprecise uncertainty measures, what degree of consistency is required for this?
Should the measure be coherent, or does some weaker concept like $2$-coherence already suffice?

In this paper, we investigate how some fundamental probability inequalities generalise in an imprecise framework.
We shall see that:
\begin{itemize}
\item[(a)]
More inequalities may correspond to a single probability inequality.
\item[(b)]
The generalised inequalities do not require exclusively coherence, even though coherence may ensure a larger applicability or tighter bounds.
\end{itemize}
These issues have been very little investigated in the literature on imprecise probabilities.
Some papers \cite{Cozman_2010,De_Cooman_Miranda_2008} introduce inequalities for sums of random numbers while studying laws of large numbers for coherent lower previsions, and in passing a Markov inequality is derived.

A different and more investigated line of research (see e.g. \cite{agahi_choquet_2018, Girotto_Holzer_2012, Mesiar_Li_Pap_2010, Roman-Flores_2007, wang_2011, zhang_2022}) extends probability inequalities to various integrals, such as the Choquet integral.
This approach partly overlaps with the lower/upper previsions one, as discussed in the later Section \ref{subsec:Jensen_without_integrals}.

The paper is organised as follows.
Section \ref{sec:preliminaries} recalls preliminary notions from the theory of Imprecise Probabilities.
In Section \ref{sec:Jensen_related_inequalities} we first derive versions of Jensen's inequality (Theorem \ref{thm:Jensen_2_coh}) for $2$-coherent previsions from an extension of a result in \cite{agahi_choquet_2018}.
Then, we generalise a recently proposed \cite{bethmann_improvement_2018} improvement of Jensen's inequality in Section \ref{subsec:improving_Jensen}.
Our result (Theorem \ref{thm:Jensen_improvement}), while assuming at least $2$-coherence, is more general than that in \cite{bethmann_improvement_2018} even in the special precise probability environment.
In Section \ref{sec:Applying Jensen's Inequalities} some applications of Jensen's imprecise inequalities are discussed.
Lyapunov-like inequalities are obtained for $2$-coherent previsions in Section \ref{subsec:Lyapunov}, while Section \ref{subsec:Basic inferences Jensen} considers some basic inferences with Jensen's inequalities, including the even moment problem.
We point out that if simple gambles only are involved and the imprecise previsions are coherent, a linear programming problem may be an attractive alternative to Jensen's inequality.
In Section \ref{sec:Markov} we obtain two Markov's inequalities and comment on their extent according to whether $2$-coherence or coherence is assumed.

In classical probability, given $X\geq 0$, Markov's inequality bounds the probability of event $(X\geq c)$, using only the expectation $E(X)$. Not requiring non-negativity for $X$ and knowing additionally the variance $\sigma^2(X)$, events like $(X\leq c)$ or $(X\geq c)$ may be bounded by Cantelli's inequalities. These are known \cite{ghosh_probability_2002} to be the sharpest ones in the family of all distributions of $X$ with given $E(X)$, $\sigma^2(X)$.
In Section \ref{sec:cantelli_inequalities} we investigate how Cantelli's inequalities can be extended to imprecise evaluations.
Firstly, we derive them in Section~\ref{subsec:cantelli_inequalities_precise} in a precise framework, but referring to linear or dF-coherent previsions rather than (less generally) expectations.
Then we obtain more forms of Cantelli's inequalities, under weak consistency assumptions, in Section \ref{subsec:cantelli_inequalities_imprecise} (Proposition \ref{pro:cantelli_imprecise}, Corollary \ref{cor:cantelli_general}). In Section \ref{subsec:cantelli_coherence} we require coherence, which lets us derive stricter Cantelli's inequalities (Proposition \ref{pro:cantelli_coherence}).
Interestingly, they involve lower and upper variances, introduced in \cite[Appendix~G]{walley_statistical_1991}.
When both can be applied, we perform a comparison between Cantelli's and Markov's inequalities, showing that neither is uniformly preferable in an imprecise framework.
Section~\ref{sec:conclusions} contains our conclusions and ideas for further developments.

\section{Preliminaries}
\label{sec:preliminaries}
Throughout the paper, we shall be concerned with bounded random variables, termed \emph{gambles}, defined on a common \emph{partition} $\prt$ of the sure event $\Omega$.
Partition $\prt$ (also called possibility space or universe of discourse) is made up of an arbitrary and possibly infinite number of exhaustive and pairwise disjoint (non-impossible) events.
A gamble $X$ is (identified by) a map from $\prt$ into $\rset$ and the \emph{image set} of $X$ is $Im(X)=\{X(\omega): \omega\in\prt\}(\subset\rset)$.
If $f$ is a real function whose domain includes $Im(X)$, $f(X)$ is a gamble with image set $Im(f(X))=\{f(X(\omega)): \omega\in\prt\}$.

We term $\asetpa$ the set of events logically dependent on $\prt$ (the powerset of $\prt$ in set-theoretical language);
an event $B$ belongs to $\asetpa$ if and only if $B$ is the logical sum of the events $\omega\in\prt$ implying $B$ ($B=\bigvee_{\omega\Rightarrow B}\omega$).
For a given $B\in\asetpa$, $B\neq\emptyset$, the conditional gamble $X|B$ is identified by the restriction of $X$ on $B$.
Thus, the image set of $X|B$ is $\{X(\omega): \omega\in\prt,\, \omega\Rightarrow B\}$, while $X|B$ is undefined for $\omega\not\Rightarrow B$.
Likewise, the image set of $f(X|B)$ is given by $Im(f(X|B))=\{f(X(\omega)): \omega\in\prt,\, \omega\Rightarrow B\}$.
\subsection{Precise Previsions}
\label{subsec:precise_previsions}
Lower and upper previsions \cite{troffaes_lower_2014,walley_statistical_1991} originate from the notion of (precise) prevision, a key concept in de Finetti's subjective probability approach \cite{de_Finetti_1970_book}.
Basically, given an arbitrary set of gambles $\dset$, a (precise) \emph{prevision} $P:\dset\rightarrow\rset$ is a mapping that associates a real number $P(X)$ that `synthesises' $X$ to each $X\in\dset$ .
To avoid inconsistent syntheses, de Finetti developed the concept of coherent prevision, which we term here \emph{dF-coherent prevision} to prevent confusion with the same term referred to lower or upper previsions:
\begin{definition}
\label{def:dF_coherent_prevision}
A mapping $P:\dset\rightarrow\rset$ is a \emph{dF-coherent prevision} iff,
$\forall n\in\natset,\forall s_0,\ldots,s_n\in\rset,\forall X_0,\ldots,X_n\in\dset$, defining $G=\sum_{i=0}^{n}s_i(X_i-P(X_i))$
we have that $\sup G\geq 0$.
\end{definition}
While this definition has a well known betting interpretation \cite{de_Finetti_1970_book,walley_statistical_1991}, in the special case that $\dset$ is a linear space $\lset$ a functional $P$ is a dF-coherent (or linear) prevision on $\lset$ iff it satisfies the following axioms \cite[Section~2.3.6]{walley_statistical_1991}:
\begin{itemize}
\item[(L)] $P(X+Y)=P(X)+P(Y), \forall X,Y\in\lset$ (linearity)
\item[(I)] $\inf X\leq P(X)\leq\sup X, \forall X\in\lset$ (internality).
\end{itemize}
These axioms imply
\begin{itemize}
\item[(H)] $P(\lambda X)=\lambda P(X), \forall\lambda\in\rset, X\in\lset$ (homogeneity)
\end{itemize}
and if $\dset$ is not a linear space, they are necessary (but not sufficient) conditions for dF-coherence of $P$.

Moreover, if a finitely additive probability $P_0$ is given on $\asetpa$, its extension to a dF-coherent prevision for $X$ is unique and coincides with the (linear) expectation of $X$ \cite[Sections 3.2.1, 3.2.2]{walley_statistical_1991}.
Yet, an agent might assess a dF-coherent prevision for any $X$ without knowing or knowing only partially either $P_0$ or even the distribution function of $X$, and in this sense the notion of dF-coherent prevision is more general than that of expectation.

\subsection{Imprecise Previsions}
\label{subsec:imprecise_previsions}
When assessing a precise prevision turns out to be difficult or unreliable, an agent may replace $P(X)$ with a lower prevision $\lpr(X)$ or an upper prevision $\upr(X)$ (or both).
In practice, we shall assume, as customary, that lower and upper previsions are \emph{conjugate}, meaning that
\begin{eqnarray}
\label{eq:conjugacy}
\upr(X)=-\lpr(-X).
\end{eqnarray}
Thus, for instance, assessing $\lpr(X)$, $\upr(X^2)$ is equivalent to assessing $\lpr$ on $\{X,-X^2\}$.
In general, conjugacy lets us refer to one type of imprecise prevision only.

The most widespread notion of consistency for lower/upper previsions is \emph{coherence}.
It is defined by introducing constraints on the sign of $s_0,\ldots,s_n$ in Definition~\ref{def:dF_coherent_prevision} \cite[Section~2.5]{walley_statistical_1991}, or by relaxing axioms (L), (I) and (H) if $\dset$ is a linear space $\lset$.
In this latter case, $\lpr$ is a coherent lower prevision on $\lset$ iff it satisfies
\cite[Section~2.5.5]{walley_statistical_1991}
\begin{itemize}
	\item[(S)] $\lpr(X+Y)\geq\lpr(X)+\lpr(Y), \forall X,Y\in\lset$ (superlinearity).
	\item[(LI)] $\inf X \leq\lpr(X), \forall X \in \lset$ (lower internality).
	\item[(PH)] $\lpr(\lambda X)=\lambda\lpr(X), \forall X\in\lset, \forall\lambda\geq 0$ (positive homogeneity).
\end{itemize}
For what follows, it is useful to recall the characterisation of coherence given by the following
\begin{theorem}[Lower Envelope Theorem]
\label{thm:lower_envelope}
$\lpr:\dset\rightarrow\rset$ is a \emph{coherent lower prevision} on $\dset$ iff there exists a (non-empty) set $\mset^*$ of dF-coherent previsions on $\dset$ such that
\begin{eqnarray}
\label{eq:envelope_inf}
\lpr(X)=\inf_{P\in\mset^*} P(X), \forall X\in\dset.
\end{eqnarray}
The infimum is attained in \eqref{eq:envelope_inf} when $\mset^*$ is the set $\mset=\{P: P\text{ is dF-coherent }\\ \text{on }\dset, P\geq\lpr\}$, termed \emph{credal set} of $\lpr$.
\end{theorem}

The Lower Envelope Theorem gives a robustness interpretation of coherent lower previsions: $\mset^*$ may be viewed as a set of potential previsions, among which the agent cannot easily establish which is the `true' one.
The credal set $\mset$ is the largest such set and is weak*-compact and convex \cite[Section 3.6.1]{walley_statistical_1991}.

A dF-coherent prevision $P$ is a special coherent lower and upper prevision such that $\lpr=\upr=P$.

A lower prevision $\lpr$ \emph{avoids sure loss} (ASL) \cite[Section 3.3.3 (a)]{walley_statistical_1991} on $\dset$ iff the credal set $\mset$ is non-empty, i.e. iff there exists a dF-coherent prevision $P\geq\lpr$ on $\dset$.
Thus, the notion of avoiding sure loss is weaker than coherence.

$2$-coherence is obtained introducing constraints in Definition~\ref{def:dF_coherent_prevision}.
Next to the additional coherence constraints on $s_0,\ldots,s_n$, it requires that $n=1$, so that $G$ is replaced by a summation of two terms.
Hence, coherence implies $2$-coherence, which is a weaker consistency requirement (unless $D$ is made of one or two gambles only, in which case they are equivalent).
We report the relevant definitions \cite{Pelessoni_Vicig_2016}:
\begin{definition}
	\label{def:2_coherent}
	\,
	\begin{itemize}
		\item[(a)]
		$\lpr:\dset\rightarrow\rset$ is a \emph{$2$-coherent lower prevision} on $\dset$ iff $\,\forall X_0,X_1\in\dset, \forall s_1\geq 0, \forall s_0\in\rset$, defining $\underline{G}_2=s_1(X_1-\lpr(X_1))-s_0(X_0-\lpr(X_0))$, we have that $\sup\underline{G}_2\geq 0$.
		\item[(b)]
		$\upr:\dset\rightarrow\rset$ is a \emph{$2$-coherent upper prevision} on $\dset$ iff $\,\forall X_0,X_1\in\dset, \forall s_1\geq 0, \forall s_0\in\rset$, defining $\overline{G}_2=s_1(\upr(X_1)-X_1)-s_0(\upr(X_0)-X_0)$, we have that $\sup\overline{G}_2\geq 0$.
	\end{itemize}
\end{definition}
$2$-coherence too can be defined by means of axioms on structured sets.
Since a limited part of the sequel will involve conditional gambles, we shall state a result in the case that $\dset$ is now the linear space $\lset|B$ of all conditional gambles $X|B$, where $X\in\lset$, a given linear space, and $B\neq\emptyset$, $B\in\asetpa$, is a fixed event.

Preliminarily, consider the following axioms, referred for later usage to a generic uncertainty measure $\mu$:
\begin{itemize}
	\setlength{\itemindent}{0.9cm}
	\item[(M)] If $X|B\geq Y|B$ then $\mu(X|B)\geq\mu(Y|B)$ (monotonicity).
	\item[(T)] $\mu(X+a|B)=\mu(X|B)+a, \forall a\in\rset$ (translation invariance).
	\item[(PH)] $\forall\lambda\geq 0, \mu(\lambda X|B)=\lambda\mu(X|B)$ (positive homogeneity).
	\item[(NH)] $\forall\lambda<0, \mu(\lambda X|B)\leq\lambda\mu(X|B)$ (negative homogeneity).	
\end{itemize}
Then we have, as a special case of \cite[Proposition 9]{Pelessoni_Vicig_2016} (cf. also \cite[Remark~1]{Pelessoni_Vicig_2016})
\begin{theorem}
\label{thm:axioms_2_coherent}
$\lpr:\lset|B\rightarrow\rset$ is a $2$-coherent lower prevision on $\lset|B$ iff axioms (M), (T), (PH), (NH) above obtain (with $\mu=\lpr$, and for all $X|B, Y|B\in\lset|B$).
\end{theorem}
\begin{remark}
\label{rem:upper_dominates_lower}
It is easy to verify that axioms (M), (T), (PH) imply further properties of $\lpr$ and its conjugate $\upr$ (cf. also \cite{Pelessoni_Vicig_2016}).
These include internality (I) (with $X$ replaced by $X|B$) and $\lpr(c|B)=\upr(c|B)=c, \forall c\in\rset$.

Importantly, it follows instead from (NH) and the conjugacy equality $\upr(X|B)=-\lpr(-X|B)$ that
\begin{eqnarray}
\label{eq:upper_dominates_lower}
\lpr(X|B)\leq\upr(X|B).
\end{eqnarray}
This seemingly obvious property is nevertheless not always implied by weaker consistency notions than $2$-coherence. 
\end{remark}
Axioms (M), (T), (PH), (NH) are necessary conditions for $2$-coherence and therefore also for coherence on a generic domain $\dset$.
Theorem \ref{thm:axioms_2_coherent} applies in particular when $B=\Omega$, i.e. when $\lpr$ is an unconditional measure on the linear space $\lset|\Omega=\lset$.

A fundamental result for $2$-coherent, coherent and dF-coherent previsions defined on a set $\dset$ is that  they allow an extension of the same type ($2$-coherent, coherent or dF-coherent, respectively) on any superset of gambles $\dset^\prime\supset\dset$ \cite{de_Finetti_1970_book,Pelessoni_Vicig_2017,walley_statistical_1991}.
Therefore, assessing one such prevision on an arbitrary set of gambles is not penalising: we can always extend it preserving its degree of consistency.

The extension is usually not unique, but a lower prevision $\lpr$ coherent on $\dset$ always has a least-committal coherent extension $\lpe$ on $\dset^\prime$, termed \emph{natural extension}.
This means that $\lpe\leq\lpr^\prime$, for any other coherent extension $\lpr^\prime$ of $\lpr$ \cite{troffaes_lower_2014,walley_statistical_1991}.
In the same manner, a $2$-coherent lower prevision $\lpr$ always admits a $2$-coherent natural extension $\lpe_2$ on any $\dset^\prime$.
Again, $\lpe_2$ is the least-committal extension of $\lpr$ among its $2$-coherent ones \cite{Pelessoni_Vicig_2016,Pelessoni_Vicig_2017}. We shall encounter natural extensions and $2$-coherent natural extensions in Section~\ref{subsec:Basic inferences Jensen}. 

A special situation occurs when a gamble is the \emph{indicator} $I_A$ of an event $A$. The prevision of an indicator is the same as the probability of the corresponding event, and these terms will be used interchangeably in this paper.
If $\dset$ is made of (indicators of) events only, then $\lpr$ ($\upr$) is a lower (upper) probability on $\dset$, whilst $P$ is a dF-coherent probability.
In this case, the conjugacy equality \eqref{eq:conjugacy} is written $\upr(A)=1-\lpr(\nega{A})$.

If the set of events $\dset$ is an algebra, a $2$-coherent $\lpr$ (or its conjugate $\upr$) is a (normalised) \emph{capacity} or \emph{fuzzy measure}, with the additional property \eqref{eq:upper_dominates_lower}.
Coherent lower/upper previsions include a number of models as special cases, such as belief functions, possibility measures \cite{troffaes_lower_2014}, coherent risk measures \cite{pelessoni_imprecise_2003}, several neighbourhood models \cite{corsato_nearly-linear_2019, montes_neighbourhood_I, montes_neighbourhood_II}, and many others.
For further information on the material recalled in this section, see among others \cite{Pelessoni_Vicig_2016,Pelessoni_Vicig_2017,troffaes_lower_2014,walley_statistical_1991}. 

\section{Jensen's and Related Inequalities}
\label{sec:Jensen_related_inequalities}
Various papers discuss extensions of Jensen's inequality beyond classical probabilities, among them \cite{Girotto_Holzer_2012, Mesiar_Li_Pap_2010,zhang_2022}. Often, they make use of Choquet, Sugeno or other integrals.
The work in \cite{agahi_choquet_2018} is the closest to our approach to Jensen's inequality with imprecise previsions, since it requires very general conditions. We start with a modified and conditional version of \cite[Theorem~2]{agahi_choquet_2018}.
For this, let $f'_{+}(x)$, $f'_{-}(x)$ be the right and left derivatives, respectively, of a real function $f$ at $x$.
It is well known that if $f$ is convex on an interval then we have, at any two interior points $x<y$ of the interval:
\begin{eqnarray}
	\label{eq:convex_property_1}
	-\infty<f'_{-}(x)\leq f'_{+}(x)\leq f'_{-}(y)\leq f'_{+}(y)<+\infty;
\end{eqnarray}
the inequalities in \eqref{eq:convex_property_1} are reversed when $f$ is concave (see e.g. \cite[Theorem~1.4.2]{Niculescu_2018}).

\begin{theorem}
\label{thm:Jensen_base}
Let $X|B$ be a gamble, $I\subset\rset$ an interval that contains the image set of $X|B$, $\phi:I\rightarrow\rset$ a convex function, $\psi:I\rightarrow\rset$ a concave function, $X|B\in\dset$.\footnote{
Here (and similarly in the next results) one may assume that $\dset$ contains all the gambles $\mu$ is applied to in the proof.
Alternatively, $\dset=\{X|B\}$ with the proviso that $\mu$ is extended to any gamble in the proof according to axioms (M), (T), (PH).
}
Further, let $\mu:\dset\rightarrow\rset$ be an uncertainty measure satisfying axioms $(M)$, $(T)$, $(PH)$ and $\nega{\mu}$ its conjugate (i.e. $\nega{\mu}(Y|B)=-\mu(-Y|B)$).
Let also $\mu(X|B)$ be an interior point of $I$.
\begin{itemize}
	\item[(a)]
	Let $J_{\phi}=[\phi^{'}_{-}(\mu(X|B)), \phi^{'}_{+}(\mu(X|B))]$.
	\begin{align}
		\label{eq:phi_1}
		\text{If } &\exists\lambda\in J_{\phi}, \lambda\geq 0 \text{ then } \mu(\phi(X|B))\geq\phi(\mu(X|B));\\
		\label{eq:phi_2}
		\text{If } &\exists\lambda\in J_{\phi}, \lambda\leq 0 \text{ then } \nega{\mu}(\phi(X|B))\geq\phi(\mu(X|B)).
	\end{align}
\item[(b)]
Let $J_{\psi}=[\psi^{'}_{+}(\mu(X|B)), \psi^{'}_{-}(\mu(X|B))]$.
	\begin{align}
		\label{eq:psi_1}
		\text{If } &\exists\lambda\in J_{\psi},\lambda\geq 0 \text{ then } \mu(\psi(X|B))\leq\psi(\mu(X|B));\\
		\label{eq:psi_2}
		\text{If } &\exists\lambda\in J_{\psi},\lambda\leq 0 \text{ then } \nega{\mu}(\psi(X|B))\leq\psi(\mu(X|B)).
	\end{align}
\end{itemize}
\end{theorem}
\begin{proof}
The proofs of $(a)$ and $(b)$ are similar. Let us prove $(b)$.
Recall for this that, since $\psi$ is concave, it holds, for any $x_0$ interior point of $I$, that
\begin{equation}
	\label{eq:convex_prop_1}
	\psi(x)\leq\psi(x_0)+\lambda (x-x_0), \forall\lambda\in [\psi^{'}_{+}(x_0), \psi^{'}_{-}(x_0)].
\end{equation}
When $x_0=\mu(X|B)$, $x=X|B$, Equation \eqref{eq:convex_prop_1} becomes
\begin{equation}
	\label{eq:convex_prop_1_rn}
	\psi(X|B)\leq\psi(\mu(X|B))+\lambda (X|B-\mu(X|B)).
\end{equation} 
\begin{itemize}
	\item[$\bullet$]
	If $\lambda\geq 0$, by axioms $(M)$, $(T)$, $(PH)$, we get from \eqref{eq:convex_prop_1_rn}
	\begin{equation*}
		\mu(\psi(X|B))\leq\psi(\mu(X|B))+\lambda (\mu(X|B)-\mu(X|B))=\psi(\mu(X|B)).
	\end{equation*}
	\item[$\bullet$]
	If $\lambda\leq 0$, rewrite \eqref{eq:convex_prop_1_rn} as
	\begin{align*}
		-\lambda (X|B-\mu(X|B))\leq\psi(\mu(X|B))-\psi(X|B).
	\end{align*}
	Now apply $(M)$, $(T)$, $(PH)$:
	\begin{align*}
		\mu(-\lambda(X|B-\mu(X|B)))&=-\lambda(\mu(X|B)-\mu(X|B))=0\\
		&\leq\psi(\mu(X|B))+\mu(-\psi(X|B)).
	\end{align*}
	Hence,
	\begin{equation*}
		-\mu(-\psi(X|B))=\nega{\mu}(\psi(X|B))\leq\psi(\mu(X|B)).
	\end{equation*}
\end{itemize}
\end{proof}
Now, a key point for the sequel is that Theorem~\ref{thm:Jensen_base} can be applied to $2$-coherent lower and upper previsions, since axioms $(M)$, $(T)$, $(PH)$  obtain for them, cf. Theorem~\ref{thm:axioms_2_coherent}.
Recalling the same Theorem~\ref{thm:axioms_2_coherent} and Remark~\ref{rem:upper_dominates_lower}, it appears clearly that $2$-coherence is the weakest consistency notion for imprecise previsions that Theorem \ref{thm:Jensen_base} can be reasonably applied to.
In fact, $2$-coherence requires, in addition to the axioms in Theorem \ref{thm:Jensen_base}, only axiom (NH), necessary to guarantee the very desirable property \eqref{eq:upper_dominates_lower}.

In classical probability theory, the commonest versions of Jensen's inequality refer to an unconditional environment.
We too shall derive now their unconditional imprecise counterparts.
For this, put $B=\Omega$ in Theorem~\ref{thm:Jensen_base}, so that $X|B=X|\Omega=X$.
\begin{theorem}[Jensen's inequalities]
\label{thm:Jensen_2_coh}
In the assumptions of Theorem \ref{thm:Jensen_base}, let further $B=\Omega$, $\{\mu, \nega{\mu}\}=\{\lpr,\upr\}$, with $\lpr$, $\upr$ $2$-coherent on their domain. Then,
\begin{align}
	\label{eq:Jensen_inf_2}
	\lpr(\psi(X))&\leq\min\{\psi(\lpr(X)),\psi(\upr(X))\},\\
	\label{eq:Jensen_inf_4}
	\upr(\phi(X))&\geq\max\{\phi(\lpr(X)),\phi(\upr(X))\}.
\end{align}	
Besides,
\begin{align}
	\begin{split}
	\label{eq:Jensen_inf_1_fork}
	\text{if }\phi'_+(\lpr(X))&\geq 0\text{ then }\lpr(\phi(X))\geq\phi(\lpr(X)),\\
	\text{if }\phi'_-(\upr(X))&\leq 0\text{ then }\lpr(\phi(X))\geq\phi(\upr(X));
	\end{split}\\
	\begin{split}
	\label{eq:Jensen_inf_3_fork}
	\text{if }\psi'_-(\upr(X))&\geq 0\text{ then }\upr(\psi(X))\leq\psi(\upr(X)),\\
	\text{if }\psi'_+(\lpr(X))&\leq 0\text{ then }\upr(\psi(X))\leq\psi(\lpr(X)).
	\end{split}
\end{align}
\end{theorem}
\begin{proof}
To prove \eqref{eq:Jensen_inf_2}, put $\mu=\lpr$ in \eqref{eq:psi_1} and \eqref{eq:psi_2} getting, respectively
\begin{align*}
	\lpr(\psi(X))\leq\psi(\lpr(X)), \text{if } \lambda\geq 0,\\
	\lpr(\psi(X))\leq\upr(\psi(X))\leq\psi(\lpr(X)), \text{if } \lambda\leq 0.
\end{align*}
Hence, no matter which is the sign of $\lambda$, it holds that
\begin{equation}
	\label{eq:to_Jensen_psi_1}
	\lpr(\psi(X))\leq\psi(\lpr(X)).
\end{equation}
Putting instead $\nega{\mu}=\lpr$ in \eqref{eq:psi_1} and \eqref{eq:psi_2} we come to the inequalities, respectively
\begin{align*}
	\lpr(\psi(X))\leq\upr(\psi(X))\leq\psi(\upr(X)), \text{if } \lambda\geq 0,\\
	\lpr(\psi(X))\leq\psi(\upr(X)), \text{if } \lambda\leq 0,
\end{align*}
and hence we conclude that
\begin{equation}
	\label{eq:to_Jensen_psi_2}
	\lpr(\psi(X))\leq\psi(\upr(X)).
\end{equation}
Grouping \eqref{eq:to_Jensen_psi_1} and \eqref{eq:to_Jensen_psi_2} we obtain \eqref{eq:Jensen_inf_2}.

We prove now \eqref{eq:Jensen_inf_1_fork}.
Put $\mu=\lpr$ in \eqref{eq:phi_1} and note that condition $\phi'_+(\lpr(X))\geq 0$ is equivalent (cf. \eqref{eq:convex_property_1}) to requiring that there is $\lambda\in J_\phi, \lambda\geq 0$.
Then, \eqref{eq:phi_1} implies $\lpr(\phi(X))\geq\phi(\lpr(X))$.

Similarly, when $\phi'_-(\upr(X))\leq 0$, apply \eqref{eq:phi_2} with $\mu=\upr$ to obtain
$\lpr(\phi(X))\geq\phi(\upr(X))$.

The proofs of \eqref{eq:Jensen_inf_4} and \eqref{eq:Jensen_inf_3_fork} are analogous to those of \eqref{eq:Jensen_inf_2} and \eqref{eq:Jensen_inf_1_fork}, respectively.
\end{proof}
\begin{remark}
\label{rem:lambda_independence}
Theorem \ref{thm:Jensen_2_coh} improves over Theorem \ref{thm:Jensen_base} in that its inequalities \eqref{eq:Jensen_inf_2}, \eqref{eq:Jensen_inf_4}  are \emph{independent} of $\lambda$.
With different consistency notions this kind of independence is no longer ensured.

It is important to observe that an additional \emph{monotonicity hypothesis} of either $\phi(x)$ or $\psi(x)$ on $[\lpr(X),\upr(X)]$ simplifies and enhances Jensen's inequalities.
Monotonicity of $\phi$ solves easily the maximum in \eqref{eq:Jensen_inf_4} and ensures that one of the equations \eqref{eq:Jensen_inf_1_fork} always applies. Similarly with $\psi$. 
\end{remark}
Another notable simplification occurs when $\lpr=\upr=P$, and $P$ is a dF-coherent prevision:
\begin{corollary}
	\label{cor:Jensen_corollary}
	Given $X$, $I$, $\phi$, $\psi$ as in Theorem \ref{thm:Jensen_base}, $B=\Omega$, let $P:\dset\rightarrow\rset$ be a dF-coherent prevision.
	Then,
	\begin{equation}
		\label{eq:corollary_jensen}
		P(\phi(X))\geq\phi(P(X)),\ P(\psi(X))\leq\psi(P(X)).
	\end{equation}
\end{corollary}
\begin{proof}
	Since $P$ is also $2$-coherent, Theorem \ref{thm:Jensen_2_coh} applies.
	Substituting $\lpr=\upr=P$ in \eqref{eq:Jensen_inf_2} and \eqref{eq:Jensen_inf_4} we get Equation \eqref{eq:corollary_jensen}.
\end{proof}
Recall that, for a gamble $X$, its expectation $E(X)$ is a coherent prevision.
Therefore, we may replace $P$ with $E$ in \eqref{eq:corollary_jensen} to obtain the familiar Jensen's inequality for convex or concave functions of $X$:
\begin{equation}
	\label{eq:traditional_Jensen}
	E(\phi(X))\geq\phi(E(X)),\ E(\psi(X))\leq\psi(E(X)).
\end{equation}
In this paper, we mostly focus on inequalities in an unconditional setting.
Yet, Theorem~\ref{thm:Jensen_base} can be exploited to derive a conditional version of Theorem~\ref{thm:Jensen_2_coh} and Corollary~\ref{cor:Jensen_corollary}, with simple modifications in their proofs.
We restrict ourselves to stating the conditional Corollary~\ref{cor:Jensen_corollary}, since it will be needed in the proof of the later Lemma \ref{lem:cantelli_preliminary}.
\begin{corollary}
	\label{cor:Jensen_corollary_cond}
	Given $X|B$, $I$, $\phi$, $\psi$ as in Theorem \ref{thm:Jensen_base}, let $P:\dset\rightarrow\rset$ be a coherent prevision.
	Then,
	\begin{equation}
		\label{eq:corollary_jensen_cond}
		P(\phi(X|B))\geq\phi(P(X|B)),\ P(\psi(X|B))\leq\psi(P(X|B)).
	\end{equation}
\end{corollary} 

\subsection{Jensen's inequalities with and without integrals}
	\label{subsec:Jensen_without_integrals}
	In classical probability theory, the expectations in Jensen's inequalities \eqref{eq:traditional_Jensen} are written as (Lebesgue) integrals.
	Quite naturally then, various papers, including \cite{Girotto_Holzer_2012,Mesiar_Li_Pap_2010,Roman-Flores_2007,wang_2011,zhang_2022}, investigated if and how Jensen's inequality can be extended to other types of integrals.
	
	As a common feature, all these approaches require an uncertainty measure $\mu$ to be preliminarily assessed on a structured subset of the powerset $\asetpa$, such as an algebra that may be $\asetpa$ itself.
	The measure $\mu$ is at least a capacity (monotone, and $\mu(\emptyset)=0$), usually normalised ($\mu(\Omega)=1$), but more restrictive conditions are sometimes assumed (for instance, lower and upper semicontinuity in \cite{Roman-Flores_2007}).
	Then, $\mu$ is extended to a set of gambles $\dset$ by means of the chosen integral.
	
	By contrast, much weaker assumptions are needed with the imprecise previsions framework of this paper:
	the measure $\mu$ may be defined \emph{only} on the set $\dset$ of gambles.
	Moreover, $\mu$ and its conjugate $\nega{\mu}$ are $2$-coherent lower and upper previsions, so that $\{\mu,\nega{\mu}\}=\{\lpr,\upr\}$.
	$2$-coherence establishes a connection with the integral approach in the following sense:
	whenever $\lpr$ or $\upr$ is ($2$-coherently) extended to $\dset\cup\asetpa$ (which is anyway not needed at all in our framework),
	the extension is a capacity on $\asetpa$. Our only extra assumption is the very reasonable condition $\lpr(\cdot)\leq\upr(\cdot)$.
	
	Studying Jensen's inequality by means of integrals brought to non-ho\-mo\-ge\-neous results.
	It is shown in \cite[Example 2.4]{Girotto_Holzer_2012} that recovering faithfully the structure of Equations \eqref{eq:traditional_Jensen} is not generally feasible with the Choquet, Sugeno or Shilkret integrals.
	With assumptions overlapping those in Theorem \ref{thm:Jensen_base}, but referring to the Choquet integral $(C)\int X d\mu=\mu(X)$,
	\cite{Girotto_Holzer_2012} obtains the inequality
	\begin{eqnarray}
		\label{GH_Jensen}
		\phi(\mu(X))=\phi\left((C) \int Xd\mu\right)\leq\max\left\{(C)\int\phi(X)d\mu,(C)\int\phi(X)d\nega{\mu}\right\}.
	\end{eqnarray}
	In \cite{zhang_2022}, a Jensen-type inequality for the Choquet integral in \cite{wang_2011} is amended by requiring monotonicity of $\phi(x)$, cf. also \cite{Mesiar_Li_Pap_2010} for analogous results under finiteness of $\prt$ and non-negativity of gambles.
	
	The classical Choquet integral satisfies axioms (M), (T), (PH) (see e.g. \cite[Appendix C]{troffaes_lower_2014}), hence these results are implied by Theorem \ref{thm:Jensen_base}.
	Equation \eqref{GH_Jensen} follows immediately from \eqref{eq:phi_1} and \eqref{eq:phi_2},
	Jensen's inequality in \cite[Theorem 2.1]{zhang_2022} from \eqref{eq:phi_1}.
	
	Other approaches, while starting from different hypotheses, are not directly comparable with ours.
	In particular, \cite{Roman-Flores_2007} obtains results with the Sugeno integral removing the convexity assumption on $\phi(x)$, while requiring its monotonicity and some further conditions.

\subsection{Improving Jensen's inequalities}
\label{subsec:improving_Jensen}
In a recent paper \cite{bethmann_improvement_2018}, an improvement to the classical Jensen's inequality \eqref{eq:traditional_Jensen} has been proposed in the case that $X$ takes values in $\zset$.
We investigate its extension to $2$-coherent imprecise previsions under more general requirements for $X$.
Roughly speaking, it suffices that the image set of $X$ has some `hole', and that $\lpr(X)$ or $\upr(X)$ `falls into' one such hole.
Let us formalise:
\begin{definition}
	\label{def:k}
	Given a gamble $X$ with image set $Im(X)$, define for $k\in[\inf X,\sup X]$
	\begin{align*}
		l(k)=\sup\{x\in Im(X):x\leq k\}\\
		u(k)=\inf\{x\in Im(X):x\geq k\}.
	\end{align*}
\end{definition}
\begin{theorem}
	\label{thm:Jensen_improvement}
	Given the gamble $X\in\dset$, the interval $I=[\inf X,\sup X]$, $\phi:I\rightarrow\rset$ convex function, $\psi:I\rightarrow\rset$ concave function, let $\lpr:\dset\rightarrow\rset$ be $2$-coherent, with $\upr$ its conjugate. Define
	\begin{equation*}
		x_L=l(\lpr(X)), x_U=u(\lpr(X)), z_L=l(\upr(X)), z_U=u(\upr(X)).
	\end{equation*}
If $\inf X\leq x_L<x_U\leq\sup X$ (in items $(a1)$, $(b2)$) and $\inf X\leq z_L<z_U\leq\sup X$ (in items $(a2)$, $(b1)$), then
\begin{enumerate}
	\item[$(a1)$]
	\begin{align}
			\label{eq:Jensen_better_1}
			\lpr(\psi(X)))&\leq\psi(x_U)\frac{\lpr(X)-x_L}{x_U-x_L}+\psi(x_L)\left(1-\frac{\lpr(X)-x_L}{x_U-x_L}\right)\\
			&\leq\psi(\lpr(X)). \notag
	\end{align}
	\item[$(a2)$]
	If $\psi(z_L)\leq\psi(z_U)$,
	\begin{align}
			\label{eq:Jensen_better_2}
			\upr(\psi(X)))&\leq\psi(z_U)\frac{\upr(X)-z_L}{z_U-z_L}+\psi(z_L)\left(1-\frac{\upr(X)-z_L}{z_U-z_L}\right)\\
			&\leq\psi(\upr(X)). \notag
	\end{align}
\item[$(b1)$]
\begin{align}
		\label{eq:Jensen_better_3}
		\upr(\phi(X)))&\geq\phi(z_U)\frac{\upr(X)-z_L}{z_U-z_L}+\phi(z_L)\left(1-\frac{\upr(X)-z_L}{z_U-z_L}\right)\\
		&\geq\phi(\upr(X)). \notag
\end{align}
\item[$(b2)$]
If $\phi(x_L)\leq\phi(x_U)$,
\begin{align}
		\label{eq:Jensen_better_4}
		\lpr(\phi(X)))&\geq\phi(x_U)\frac{\lpr(X)-x_L}{x_U-x_L}+\phi(x_L)\left(1-\frac{\lpr(X)-x_L}{x_U-x_L}\right)\\
		&\geq\phi(\lpr(X)). \notag
\end{align}
\end{enumerate}
\end{theorem}
\begin{proof}
We prove $(a1)$.
Let $y=s(x)$ be the straight line joining $(x_L,\psi(x_L))$, $(x_U,\psi(x_U))$, 
\begin{equation}
	\label{eq:str_line}
	s(x)=\frac{\psi(x_U)-\psi(x_L)}{x_U-x_L}(x-x_L)+\psi(x_L). 
\end{equation}
Since $\psi$ is concave, it holds that $\psi(x)\leq s(x)$ for both $x\leq x_L$ and $x\geq x_U$.\footnote{
This standard result can be derived for instance using the \emph{three chords inequality} for concave functions (cf. \cite[Equation (1.16)]{Niculescu_2018}).
}
On the other hand, $X$ takes no value in $]x_L,x_U[$, thus it is true that
\begin{equation}
	\label{eq:psi_vs_s}
	\psi(X)\leq s(X).
\end{equation}
From \eqref{eq:psi_vs_s}, we obtain \eqref{eq:Jensen_better_1} through the derivation that follows.
Its first inequality is due to monotonicity of $\lpr$, the second to \eqref{eq:str_line}, $(T)$ and to axiom $(PH)$ if $\psi(x_U)\geq\psi(x_L)$,
to $(NH)$
otherwise\footnote{Note that when $\psi(x_U)\geq\psi(x_L)$, the second inequality is an equality.
},
the third to concavity of $\psi$:
\begin{align*}
	\lpr(\psi(X))&\leq\lpr(s(X))\leq\frac{\psi(x_U)-\psi(x_L)}{x_U-x_L}(\lpr(X)-x_L)+\psi(x_L)\\
	&=\psi(x_U)\frac{\lpr(X)-x_L}{x_U-x_L}+\psi(x_L)\left(1-\frac{\lpr(X)-x_L}{x_U-x_L}\right)\leq\psi(\lpr(X)).
\end{align*}
To prove $(a2)$, follow the same course of reasoning as in $(a1)$,
replacing $x_L, x_U, \lpr$ with $z_L, z_U, \upr$ and applying $(PH)$ at the second inequality\footnote{
This is the only asymmetry with the proof of $(a1)$: when $\psi(z_L)>\psi(z_U)$, we cannot apply the upper prevision version of $(NH)$ $\upr(\lambda X)\geq\lambda\upr(X)$, because the inequality direction is opposite to that of the majorisation we would need.
}
of the final derivation.

The proof of $(b1)$, $(b2)$ is analogous to $(a1)$, $(a2)$ respectively.
\end{proof}

How can Theorem~\ref{thm:Jensen_improvement} be exploited, operationally?
To exemplify, let us discuss $(a1)$ and $(a2)$.
Note that the hypothesis $\inf X\leq x_L< x_U\leq\sup X$
\begin{itemize}
	\item[(a)]
	is compatible with the assessments $\lpr(X)=x_L>\inf X$ or $\lpr(X)=x_U<\sup X$.
\end{itemize}
Besides,
\begin{itemize}
 	\item[(b)]
 	the values $\lpr(X)=x_L$ or $\lpr(X)=x_U$ make \eqref{eq:Jensen_better_2} useless.
	In fact, in this case, \eqref{eq:Jensen_better_1} boils down 
	to Equation~\eqref{eq:to_Jensen_psi_1}, $\lpr(\psi(X))\leq\psi(\lpr(X))$, and cannot improve any of Jensen's inequalities.
\end{itemize}
Condition $\inf X\leq z_L<z_U\leq\sup X$ plays the same role with $\upr$ in $(a2)$ as for $(a)$, but not for $(b)$ as far as Equations \eqref{eq:Jensen_inf_3_fork} are concerned, as will appear shortly.

Items $(a1)$ and $(a2)$ may tighten Jensen's inequalities \eqref{eq:Jensen_inf_2} and \eqref{eq:Jensen_inf_3_fork}.
To see this, let us label the right-hand terms of the inequalities \eqref{eq:Jensen_better_1} and \eqref{eq:Jensen_better_2}:
\begin{align*}
	M_1&=\psi(x_U)\frac{\lpr(X)-x_L}{x_U-x_L}+\psi(x_L)\left(1-\frac{\lpr(X)-x_L}{x_U-x_L}\right)\\
	M_2&=\psi(z_U)\frac{\upr(X)-z_L}{z_U-z_L}+\psi(z_L)\left(1-\frac{\upr(X)-z_L}{z_U-z_L}\right)
\end{align*}
Consider inequality \eqref{eq:Jensen_inf_2}.
By $(a1)$, we may replace it with
\begin{align}
	\label{eq:M1_is_better}
	\lpr(\psi(X))\leq\min\{M_1, \psi(\upr(X))\},	
\end{align} 
since $M_1\leq\psi(\lpr(X))$.
Note that if $\psi$ is monotone non-decreasing on $[\lpr(X),\upr(X)]$, or just when $\psi(\lpr(X))\leq \psi(\upr(X))$,
\eqref{eq:M1_is_better} simplifies to $\lpr(\psi(X))\leq M_1$.

When $(a2)$ is applicable, we have that
\begin{align}
	\label{eq:M2_is_better}
	\lpr(\psi(X))\leq\upr(\psi(X))\leq M_2\leq\psi(\upr(X)).
\end{align}
This brings to the next tightening of \eqref{eq:M1_is_better} and hence of inequality \eqref{eq:Jensen_inf_2}:
\begin{align*}
	\lpr(\psi(X))\leq\min\{M_1,M_2\}.
\end{align*}
Now take inequalities \eqref{eq:Jensen_inf_3_fork}.
Here the extent of Theorem \ref{thm:Jensen_improvement} is more limited:
\emph{if} $(a2)$ can be applied, from \eqref{eq:M2_is_better}, the first inequality \eqref{eq:Jensen_inf_3_fork}
can be replaced by $\upr(\psi(X))\leq M_2$, i.e. by \eqref{eq:Jensen_better_2}.
 
Note that \eqref{eq:Jensen_better_2} may apply, while the first inequality \eqref{eq:Jensen_inf_3_fork} can not, when $\psi'_-(\upr(X))<0$.
Since then $\psi'_+(\upr(X))<0$ (cf. \eqref{eq:convex_property_1}), using the second inequality \eqref{eq:Jensen_inf_3_fork} we obtain in this case that
\begin{eqnarray*}
	\upr(\psi(X))\leq\min\{M_2,\psi(\lpr(X))\},
\end{eqnarray*}
structurally similar to \eqref{eq:M1_is_better}.

The role of $(b1)$, $(b2)$ in improving Jensen's inequalities \eqref{eq:Jensen_inf_4}, \eqref{eq:Jensen_inf_1_fork} is analogous.

Finally, when $\lpr=\upr=P$, $P$ dF-coherent, \eqref{eq:Jensen_better_3} and \eqref{eq:Jensen_better_4} boil down to a single inequality, like  \eqref{eq:Jensen_better_1} and \eqref{eq:Jensen_better_2}. They can by applied to possibly improve Jensen's inequalities \eqref{eq:corollary_jensen}.

The special case of Theorem \ref{thm:Jensen_improvement} established in \cite{bethmann_improvement_2018} is obtained assuming additionally that an expectation $E(X)$ is given, that $Im(X)\subset\zset$, and (whenever an improvement over Jensen's inequality can be achieved) that $E(X)$ is non-integer.
In particular, in \cite{bethmann_improvement_2018} $l(E(X))$ is set equal to the floor value $\lfloor E(X)\rfloor$, $u(E(X))$ to $\lceil E(X)\rceil=\lfloor E(X)\rceil+1$, so that $x_U-x_L=z_U-z_L=1$.

\section{Applying Jensen's Inequalities}
\label{sec:Applying Jensen's Inequalities}
\subsection{Lyapunov's Inequalities}
\label{subsec:Lyapunov}
In classical probability theory, \emph{Lyapunov's inequality} ensures that (see e.g. \cite{shiryaev_probability_1996})
\begin{align*}
	\text{for }0<s<t, [E(|X|^s)]^{\frac{1}{s}}\leq [E(|X|^t)]^{\frac{1}{t}},
\end{align*}
provided that the expectations above are finite.
The result may be obtained from Jensen's inequality \cite{shiryaev_probability_1996}.
Thus, a natural question is whether Lyapunov-like inequalities can be derived in Imprecise Probability theory from (some of) Jensen's inequalities \eqref{eq:Jensen_inf_2}$\div$\eqref{eq:Jensen_inf_3_fork}.
The next proposition gives an affirmative answer.
\begin{proposition}
	\label{pro:Lyapunov}
	Let $\lpr$ be $2$-coherent wherever it is defined,
	$\upr$ being its conjugate.
	Letting $X$ be a gamble and $0<s<t$, we have that
	\begin{align}
		\label{eq:imprecise_Lyapunov_1}
		\text{(a) }&[\lpr(|X|^s)]^{\frac{1}{s}}\leq [\upr(|X|^s)]^{\frac{1}{s}}\leq [\upr(|X|^t)]^{\frac{1}{t}}.\\
		\label{eq:imprecise_Lyapunov_2}
		\text{(b) }& \text{If } X\geq 0, \text{ it holds that }
		[\lpr(X^s)]^{\frac{1}{s}}\leq[\lpr(X^t)]^{\frac{1}{t}}.
	\end{align}	
\end{proposition}
\begin{proof}\,	
\begin{itemize}
\item[(a)]
The first inequality in \eqref{eq:imprecise_Lyapunov_1} derives from $(M)$ and the property \eqref{eq:upper_dominates_lower} of $2$-coherence, ensuring that $0\leq\lpr(|X|^s)\leq\upr(|X|^s)$, after raising its terms to $\frac{1}{s}$.

To prove that
\begin{align}
	\label{eq:Lyapunov_second_part}
	[\upr(|X|^s)]^{\frac{1}{s}}\leq [\upr(|X|^t)]^{\frac{1}{t}}
\end{align}
let $r=\frac{t}{s}>1$, $\phi(x)=|x|^r$.

Noting that $\phi$ is a convex function and making use of \eqref{eq:Jensen_inf_4} at the second inequality, we obtain:
\begin{align}
	\label{eq:Lyapunov_intermediate}
	0\leq|\upr(Y)|^r=\phi(\upr(Y))\leq\upr(\phi(Y))=\upr(|Y|^r).
\end{align}
Now, putting $Y=|X|^s$ and replacing $r$ in \eqref{eq:Lyapunov_intermediate},
\begin{align*}
	|\upr(Y)|^r=(\upr(|X|^s))^\frac{t}{s}\leq\upr(|X|^{s\cdot\frac{t}{s}})=\upr(|X|^{t}),
\end{align*}
implying \eqref{eq:Lyapunov_second_part}.
\item[(b)]
Define $r$ as in $(a)$, and the convex function $\phi(x)=x^r$.
Clearly, $\phi(x)$ is increasing for $x\geq 0$.
Because of this, and since $\upr(Y)\geq\lpr(Y)\geq 0$ for any $Y\geq 0$, we have that
\begin{align*}
	(\lpr(Y))^r=\phi(\lpr(Y))\leq\lpr(\phi(Y))=\lpr(Y^r),
\end{align*}
where the inequality is due to the first of \eqref{eq:Jensen_inf_1_fork}.
The proof continues now like that of $(a)$: take $Y=X^s$ and replace $r=\frac{t}{s}$ in the last derivation to obtain
$(\lpr(X^s))^\frac{t}{s}\leq\lpr(X^t)$, and then \eqref{eq:imprecise_Lyapunov_2}.
\end{itemize}
\end{proof}
We observe that again a probability inequality is split into more imprecise probability inequalities and that the weak requirement of $2$-coherence is enough to obtain them.

As a by-product of Lyapunov's inequalities we have:
\begin{corollary}
\label{cor:Lyapunov}
If $X\geq 0$, then for $2$-coherent and conjugate $\lpr$, $\upr$ we have:
\begin{align}
\label{eq:Lyapunov_corollary}
[\lpr(X)]^2\leq\lpr(X^2), [\upr(X)]^2\leq\upr(X^2).
\end{align}
\end{corollary}
\begin{proof}
Put $s=1$, $t=2$ in Proposition \ref{pro:Lyapunov}.
This specialises \eqref{eq:imprecise_Lyapunov_2} into
$(0\leq)\leq\lpr(X)\leq[\lpr(X^2)]^\frac{1}{2}$.
Squaring gives the first inequality in \eqref{eq:Lyapunov_corollary}.
Taking instead the second inequality in \eqref{eq:imprecise_Lyapunov_1} as a starting point we obtain the remaining inequality in \eqref{eq:Lyapunov_corollary}. 
\end{proof}
Note that inequalities \eqref{eq:Lyapunov_corollary} in the form
$\mu(X^2)-[\mu(X)]^2\geq 0$ ($\mu=\lpr$ or alternatively $\mu=\upr$) appear to be an imprecise probability generalisation of the variance property $E(X^2)-E^{2}(X)=\sigma^2 (X)\geq 0$.

\subsection{Basic Inferences with Jensen's Inequalities}
\label{subsec:Basic inferences Jensen}
The most natural \emph{basic inferential problem} a Jensen's inequality can be applied to is the following:
a \emph{coherent} lower (upper) prevision $\lpr$ ($\upr$) is given on $\{X\}$ and some bound is sought for its \emph{coherent} extensions to a new gamble $Y$, which is a convex or concave function of $X$.
We point out that an agent's uncertainty information may be very essential here and limited to just $\lpr(X)$ or $\upr(X)$.
For instance, these values may have been communicated by an expert, without further details on how they were obtained.

This problem includes the imprecise version of the \emph{even moment problem}, where $Y=X^{2n}, n\in\natset^+$.
In fact, since $Y=\phi(X)$, with $\phi(x)=x^{2n}$ being a convex function, if $\lpr(X)$ and $\upr(X)$ are both assessed Equation \eqref{eq:Jensen_inf_4} returns a lower bound on $\upr(X^{2n})$.

Some further result is possible under additional hypotheses.
For instance, assessing only $\lpr(X)$ may be enough if $X\geq 0$, since then $\phi$ is increasing and the first inequality in \eqref{eq:Jensen_inf_1_fork} ensures that $\lpr(Y)\geq(\lpr(X))^{2n}$.
Still while $X\geq 0$, the \emph{odd moment problem} with $Y=X^{2n+1}$, $n\in\nset^+$, can be tackled too applying \eqref{eq:Jensen_inf_1_fork},
giving $\lpr(Y)\geq(\lpr(X))^{2n+1}$.

Returning to the general \emph{basic inferential problem}, an interesting question is:
how useful or convenient is answering it via Jensen's inequalities, even when monotonicity of $\phi$, $\psi$ or other assumptions simplify them?
The answer depends on how difficult it is to find an extension of $\lpr$ or $\upr$ to $Y$ with alternative methods.

When $X$ is a \emph{simple} gamble, hence taking finitely many values $x_1, x_2,\ldots,x_n$, we may detect the \emph{natural extension} $\lpe(Y)$, i.e. the smallest coherent extension of $\lpr$ to $Y$, via linear programming (LP).
The LP problem to solve is
\footnote{
To exemplify we consider function $\phi$. It could be replaced by $\psi$.
}
\begin{mini}
	{}{E_P(Y)=\sum_{i=1}^{n} p_{i}\phi(x_{i})}{\label{eq:dominating}}{}
	\addConstraint{ \sum_{i=1}^{n}p_i x_i=E_P(X)\geq\lpr(X)}{}{}
	\addConstraint{\sum_{i=1}^{n}p_i=1,\ p_i\geq 0\ (i=1,\dots, n)}{}{.}
\end{mini}
Its unknowns are $p_1,\ldots,p_n$, with $p_i$ probability of $(X=x_i)$.
The LP problem determines $\lpe(Y)$ as the minimum expectation $E_P(Y)$ for $Y=\phi(X)$,
among all probabilities $P$ that guarantee the dominance condition $E_P(X)\geq\lpr(X)$.
The fact that the LP solution always exists and is the natural extension $\lpe(Y)$ relies on well known properties of the natural extension (cf. \cite[Sections~3.1.1, 3.3.3, 3.4.1]{walley_statistical_1991}).
We will not detail the explanation, since this procedure is tangential in this paper.

\begin{example}
	\label{ex:improvement}
	Consider the partition $\prt=\{\omega_1,\omega_2, \omega_3\}$ and the gamble $X$ defined by $X(\omega_1)=x_1=-1$, $X(\omega_2)=x_2=1$, $X(\omega_3)=x_3=2$.
	Next to this, we only know that $\lpr(X)=0.75$.
	
	We may compute the natural extension $\lpe(X^2)$ of $\lpr$ to $\phi(X)=X^2$ solving the LP problem \eqref{eq:dominating}, i.e.	\begin{mini*}
		{}{E_P(X^2)=p_1+p_2+4 p_3}{}{}
		\addConstraint{-p_1+p_2+2 p_3\geq 0.75}{}{}
		\addConstraint{p_1+p_2+p_3=1, \ p_i\geq 0\ (i=1,2,3)}{}{.}
	\end{mini*}
	The set of optimal solutions is made of the triples $(1-p_2, p_2, 0)$, with $p_2\in [\frac{7}{8},1]$.
	Consequently, $\lpe(X^2)=1-p_2+p_2+0=1$.
	
	We may alternatively apply the first Jensen's inequality \eqref{eq:Jensen_inf_1_fork}.
	This gives the bound $\lpe(X^2)\geq(\lpr(X))^2=0.5625$.
	
	It is also possible to apply equation \eqref{eq:Jensen_better_4}, since $x_L=-1$, $x_U=1$ and $\phi(x_L)=\phi(x_U)=1$.
	We obtain the (best) bound $\lpe(X^2)\geq 1$.
\end{example}
Note that, in the LP approach we are discussing,
we need only assessing $\lpr$ on $X$ as an input, while Jensen's inequalities may require knowing also $\upr(X)$.
Moreover, after performing the LP problem we obtain a lower prevision on two gambles only, $X$ and $Y$.
Therefore coherence and $2$-coherence are not distinguishable here, and the LP problem detects also the $2$-coherent natural extension $\lpe_2$.
In fact, in this case $\lpe_2=\lpe$, as appear from their definitions in \cite[Definition 4]{Pelessoni_Vicig_2017}.

If, more generally, $\lpr$ is given on a finite set $\dset$ of simple gambles (including $X$, but not $Y$), we can still find out $\lpe(Y)$ via linear programming.
It is necessary that $\lpr$ is coherent or at least avoids sure loss on $\dset$.
Given this, it suffices to apply, in the new LP problem, the dominance constraint \eqref{eq:dominating} to any gamble in $\dset$.
However, we generally no longer obtain $\lpe_2 (Y)$ in this way, but only an upper bound.
In fact, it is known that $\lpe_2\leq\lpe$ \cite[Lemma 1]{Pelessoni_Vicig_2017}.

Clearly, the LP alternative does not apply to non-simple gambles.
If our beliefs are encoded by some specific models, (relatively) simple formulae for computing $\lpe(Y)$ may be available.
This requires anyway knowing more than just $\lpr(X)$ or $\upr(X)$.
For instance, let a $2$-monotone lower probability be assessed on the powerset $\asetpa$ of partition $\prt$.
Then, the natural extension $\lpe(Y)$ is obtained as a Choquet integral \cite{troffaes_lower_2014}.
In general, however, it is precisely with non-simple gambles that Jensen's inequalities are operationally most useful.
The improved inequalities of Theorem \ref{thm:Jensen_improvement} are instead useless if, for instance, the image of $X$ is an interval, but can be helpful in other situations (including the basic problem, with $X$ simple).

\section{Markov's Inequalities}
\label{sec:Markov}
Markov's inequality, very well known in classical probability theory, operates under rather mild requirements:
if $X$ is a non-negative random number and its expectation $E(X)$ exists, then
\begin{align*}
\forall a>0, P(X\geq a)\leq\frac{E(X)}{a}.
\end{align*}
The interesting fact we are going to point out is that Markov's inequality does not really depend on our employing precise probabilities, but is still valid with much weaker uncertainty measures such as $2$-coherent imprecise previsions.
This ensues from the next proposition.
\begin{proposition}[Markov's inequalities]
\label{pro:Markov_inequalities}
Let $X$ be a non-negative gamble, $X\geq 0$.
\begin{itemize}
\item[(a)]
Let $\lpr$ be a $2$-coherent lower prevision on $\{X, (X\geq a)\}$.
Then,
\begin{align}
\label{eq:lower_Markov_inequality}
\lpr(X\geq a)\leq\frac{\lpr(X)}{a}, \forall a>0 \text{  (Lower Markov Inequality)}.
\end{align}
\item[(b)]
Let $\upr$ be a $2$-coherent upper prevision on $\{X, (X\geq a)\}$.
Then,
\begin{align}
	\label{eq:upper_Markov_inequality}
	\upr(X\geq a)\leq\frac{\upr(X)}{a}, \forall a>0 \text{  (Upper Markov Inequality)}.
\end{align}
\end{itemize}
\end{proposition}
\begin{proof}
\begin{itemize}
\item[(a)]
Apply Definition \ref{def:2_coherent} $(a)$ with $s_1=a$, $s_0=1$, $X_1=I_{(X\geq a)}$, $X_0=X$.
We obtain
\begin{align*}
\underline{G}_2&=a(I_{(X\geq a)}-\lpr(X\geq a))-(X-\lpr(X))\\
&=\lpr(X)-a\lpr(X\geq a)+aI_{(X\geq a)}-X.
\end{align*}
Note that $aI_{(X\geq a)}-X\leq 0$.
Therefore, it is necessary for the $2$-coherence condition $\sup \underline{G}_2\geq 0$ that $\lpr(X)-a\lpr(X\geq a)\geq 0$, which is the Lower Markov Inequality.
\item[(b)]
The proof is analogous to $(a)$, using Definition \ref{def:2_coherent} $(b)$ with $s_1=1$, $s_0=a$, $X_1=X$, $X_0=I_{(X\geq a)}$.  
\end{itemize}
\end{proof}
Obviously, Markov's inequalities apply to coherent upper/lower previsions too.
Yet, deducing that working with coherence or $2$-coherence is always indifferent as far as Markov's inequalities are concerned would be erroneous.
In fact, coherent imprecise previsions satisfy a greater number of properties.
This makes the use of Markov's inequalities possible for them in more general situations than with $2$-coherence. See for this the next example.
\begin{example}
\label{ex:Markov}
Let $X,Y,Z$ be non-negative gambles such that $\upr(X)=60$, $\upr(Y)=50$, $\upr(Z)=10$.
If $\upr$ is coherent, we can get a bound to $\upr(X+Y+Z\geq 150)$ with the following derivation, making use of the upper Markov's inequality \eqref{eq:upper_Markov_inequality} firstly, and then of the sublinearity property of coherent upper previsions \cite[Section 2.6.1 (m)]{walley_statistical_1991}:
\begin{align*}
\upr(X+Y+Z\geq 150)\leq\frac{\upr(X+Y+Z)}{150}\leq\frac{\upr(X)+\upr(Y)+\upr(Z)}{150}=\frac{4}{5}.
\end{align*}
However, sublinearity is not guaranteed with $2$-coherence, which does not therefore allow the above computation.
\end{example}

\section{Cantelli's inequalities}
\label{sec:cantelli_inequalities}
Markov's inequality is probably the best known instance of one-sided inequality, bounding the probability of $(X\geq a)$.
Further known inequalities of this kind assume knowledge of some higher order moments.

In a form investigated among others by Uspensky \cite[p. 198]{uspensky_1937}, Cantelli's inequalities\footnote{
Also termed Chebyshev-Cantelli inequalities or one-sided Chebyshev inequalities, as their earliest clues appear in Chebyshev's work. 
Cantelli \cite{cantelli_1929} investigated them in 1928 (possibly also in some now hardly available paper around 1910), Uspensky \cite{uspensky_1937} in 1937.
}
apply to a random number $X$ with expectation $E(X)=0$ and variance $\sigma^2=E(X^2)$. Then, $\forall\varepsilon>0$,
\begin{eqnarray}
	\label{eq:cantelli_simple_precise_1}
	P(X\leq-\varepsilon)\leq\frac{E(X^2)}{E(X^2)+\varepsilon^2}\\
	\label{eq:cantelli_simple_precise_2}
	P(X\geq\varepsilon)\leq\frac{E(X^2)}{E(X^2)+\varepsilon^2}.
\end{eqnarray}
Inequality \eqref{eq:cantelli_simple_precise_2} is similar in its structure to Markov's.
While requiring knowledge of a second moment, it is not restricted to non-negative random numbers and is never vacuous, unlike Markov's inequality.
Inequality \eqref{eq:cantelli_simple_precise_1} is a left-sided variant.

When $E(X)\neq 0$, both \eqref{eq:cantelli_simple_precise_1} and \eqref{eq:cantelli_simple_precise_2} obtain replacing $X$ with $X-E(X)$:
\begin{eqnarray}
	\label{eq:cantelli_precise_prev_non_zero}
	P(X\leq E(X)-\varepsilon)\leq\frac{\sigma^{2}(X)}{\sigma^{2}(X)+\varepsilon^2}, 	P(X\geq E(X)+\varepsilon)\leq\frac{\sigma^{2}(X)}{\sigma^{2}(X)+\varepsilon^2}.
\end{eqnarray}

\subsection{Cantelli's inequalities with dF-coherent previsions}
\label{subsec:cantelli_inequalities_precise}
Cantelli's inequalities \eqref{eq:cantelli_simple_precise_1}, \eqref{eq:cantelli_simple_precise_2}, \eqref{eq:cantelli_precise_prev_non_zero} generalise to imprecise previsions.
Before seeing this, we take an intermediate step and prove a version of \eqref{eq:cantelli_simple_precise_1}, \eqref{eq:cantelli_simple_precise_2}, where $X$ is a gamble and $E(X)$, $E(X^2)$ are replaced by the dF-coherent previsions $P(X)$, $P(X^2)$.
A preliminary lemma is needed:\footnote{
	Our proof applies Jensen's inequality \eqref{eq:corollary_jensen_cond}. Alternatively, the Cauchy-Schwarz inequality may be employed.
}
\begin{lemma}
\label{lem:cantelli_preliminary}
Given a gamble $X$, event $A=(X>0)$ with indicator $I_A$, a dF-coherent prevision $P$ on $\{A, I_{A}X, I_{A}X^{2}\}$, it holds that
\begin{eqnarray}
	\label{eq:cantelli_preliminary}
	[P(I_{A}X)]^{2}\leq P(A)P(I_{A}X^2).
\end{eqnarray}
\end{lemma}
\begin{proof}
Suppose first $P(A)=0$. From $0\leq I_{A}X\leq I_{A}\sup X$,
monotonicity and homogeneity of $P$ imply that $0\leq P(I_{A}X)\leq P(A)\sup X=0$.
Therefore, $P(I_{A}X)=P(A)=0$ and \eqref{eq:cantelli_preliminary} holds.

Let now $P(A)>0$. We have that
\begin{eqnarray*}
[P(I_{A}X)]^2\leq P(A)P(I_{A}X^{2})\text{ iff } &[\frac{P(I_{A}X)}{P(A)}]^2\leq\frac{P(I_{A}X^2)}{P(A)}\\
\text{ iff } &[P(X|A)]^2\leq P(X^2|A), 
\end{eqnarray*}
which is a true inequality, as follows applying \eqref{eq:corollary_jensen_cond} with $\phi(x)=x^2$.
\end{proof}
\begin{proposition}
\label{pro:cantelli_precise}
Given a gamble $X$, $\varepsilon>0$, and a dF-coherent prevision $P$ on $\{X, X^2, (X\leq -\varepsilon)\}$ such that $P(X)=0$, 
\begin{align}
	\label{eq:cantelli_precise_1}
	P(X\leq-\varepsilon)\leq\frac{P(X^2)}{P(X^2)+\varepsilon^2},\\
	\label{eq:cantelli_precise_2}
	P(X\geq\varepsilon)\leq\frac{P(X^2)}{P(X^2)+\varepsilon^2}.
\end{align}
\end{proposition}
\begin{proof}
We prove first inequality \eqref{eq:cantelli_precise_1}.
Define the event $A=(X>-\varepsilon)=(X+\varepsilon>0)$.
Since $I_{A}(X+\varepsilon)\geq X+\varepsilon$, applying monotonicity and translation invariance to (any dF-coherent extension to $I_{A}(X+\varepsilon)$,  $X+\varepsilon$ of) $P$ we obtain
\begin{eqnarray}
	\label{eq:usp_proof_interm}
	P(I_{A}(X+\varepsilon))\geq P(X+\varepsilon)= \varepsilon>0.
\end{eqnarray}
By \eqref{eq:usp_proof_interm} at the first inequality, \eqref{eq:cantelli_preliminary} at the second, and noting that $I_{A}(X+\varepsilon)^2\leq (X+\varepsilon)^2$ at the third, it turns out that
\begin{align*}
	\varepsilon^2&\leq [P(I_{A}(X+\varepsilon))]^2\leq P(A)P(I_{A}(X+\varepsilon)^2)\leq P(A)P((X+\varepsilon)^2)\\
	&=P(A)(P(X^2)+\varepsilon^2)=(1-P(X\leq-\varepsilon))(P(X^2)+\varepsilon^2),
\end{align*}
using linearity and homogeneity of $P$ at the equalities.
Inequality \eqref{eq:cantelli_precise_1} follows straightforwardly.

The proof of inequality \eqref{eq:cantelli_precise_2} ensues from \eqref{eq:cantelli_precise_1} applied to $-X$ (noting that $P(-X)=-P(X)=0$ and $(X\geq\varepsilon)=(-X\leq-\varepsilon)$).
\end{proof}
We point out that Cantelli's inequalities in Proposition \ref{pro:cantelli_precise} are more general than their classical counterparts \eqref{eq:cantelli_simple_precise_1} and \eqref{eq:cantelli_simple_precise_2}.
In fact, they require no probability assessment on (subsets of) $\asetpa$.
Further, if such an assessment is available, it may be a dF-coherent probability, not necessarily countably additive. 

\subsection{Cantelli's inequalities with imprecise previsions}
\label{subsec:cantelli_inequalities_imprecise}
We focus at first on establishing an imprecise version of Proposition \ref{pro:cantelli_precise} while assuming very loose consistency requirements for the relevant imprecise previsions.
\begin{proposition}[Cantelli's inequalities]
\label{pro:cantelli_imprecise}
Let $X$ be a gamble, and take $\varepsilon>0$.
\begin{itemize}
\item[(a)]
Let $\lpr$ be defined on $\dset^{\leq}_{0}(X)=\{(X\leq-\varepsilon), -X^2\}$ and define
$\mathcal{M}_{0}^{\leq}(X)=\{P: P \text{ dF-coherent}, P\geq\lpr\text{ on }\dset^{\leq}_{0}(X)\text{ and }P(X)=0\}$.
If $\mathcal{M}_{0}^{\leq}(X)\neq\varnothing$, then
\begin{eqnarray*}
\label{cantelli_imprecise_1}
\lpr(X\leq-\varepsilon)\leq\frac{\upr(X^2)}{\upr(X^2)+\varepsilon^2}.
\end{eqnarray*}
\item[(b)]
Let $\lpr$ be defined on $\dset^{\geq}_{0}(X)=\{(X\geq\varepsilon), -X^2\}$ and define
$\mathcal{M}_{0}^{\geq}(X)=\{P: P \text{ dF-coherent}, P\geq\lpr\text{ on }\dset^{\geq}_{0}(X)\text{ and }P(X)=0\}$.
If $\mathcal{M}_{0}^{\geq}(X)\neq\varnothing$, then
\begin{eqnarray*}
	\label{cantelli_imprecise_2}
	\lpr(X\geq\varepsilon)\leq\frac{\upr(X^2)}{\upr(X^2)+\varepsilon^2}.
\end{eqnarray*}
\end{itemize}
\end{proposition}
\begin{proof}
\begin{itemize}
	\item[(a)]
	Take $P\in\mathcal{M}_{0}^{\leq}(X)$.
	The proof follows from the derivation
	\begin{align*}
		\lpr(X\leq-\varepsilon)&\leq P(X\leq-\varepsilon)\leq\frac{P(X^2)}{P(X^2)+\varepsilon^2}\\
		&\leq\frac{-\lpr(-X^2)}{-\lpr(-X^2)+\varepsilon^2}=\frac{\upr(X^2)}{\upr(X^2)+\varepsilon^2},
	\end{align*}
	using \eqref{eq:cantelli_precise_1} at the second inequality.
	For the last inequality, observe first that $P(X^2)\leq -\lpr(-X^2) \text{ iff } -P(X^2)=P(-X^2)\geq \lpr(-X^2)$,
	which obtains because $P\in\mathcal{M}_{0}^{\leq}(X)$.
	Then, apply the arithmetic inequality $\frac{a}{a+b}\leq\frac{c}{c+b}$ iff $a\leq c$, when $b(a+b)(c+b)>0$ (with $a=P(X^2)$, $c=-\lpr(-X^2)$, $b=\varepsilon^2$).
	\item[(b)]
	The derivation is similar to (a), taking $P\in\mathcal{M}_{0}^{\geq}(X)$ and using \eqref{eq:cantelli_precise_2}.
\end{itemize}
\end{proof}
We point out that $\lpr$ is asked, in Proposition \ref{pro:cantelli_imprecise} (a), to \emph{avoid sure loss} (ASL) on $\dset^{\leq}_{0}(X)$, since  $\mset_{0}^{\leq}(X)\neq\varnothing$ implies that the credal set of $\lpr$ is non-empty (cf.~Section~\ref{subsec:imprecise_previsions}).
If $\lpr$ is $2$-coherent on $\dset^{\leq}_{0}(X)$, it is also coherent and ASL on $\dset^{\leq}_{0}(X)$, but $2$-coherence is not really required in (a).
Instead, there is the additional condition $P(X)=0$ for $P\in\mset_{0}^{\leq}(X)$.
The assumptions of Proposition \ref{pro:cantelli_imprecise} (b) are analogous.

We can further generalise Cantelli's inequalities in the following way, derived straightforwardly from Proposition \ref{pro:cantelli_imprecise}.
\begin{corollary}
\label{cor:cantelli_general}
Given a gamble $X$ and $c\in\rset$, $\varepsilon>0$,
\begin{itemize}
\item[(a)]
Let $\lpr$ be defined on $\dset_{\leq}(X)=\{(X\leq c-\varepsilon), -(X-c)^2\}$, and let $\mset_{\leq}(X)=\{P: P \text{ dF-coherent}, P\geq\lpr \text{ on } \dset_{\leq}(X)\text{ and }P(X)=c\}$.
If $\mset_{\leq}(X)\neq\varnothing$, then
\begin{eqnarray}
\label{eq:cantelli_general_1}
\lpr(X\leq c-\varepsilon)\leq\frac{\upr((X-c)^2)}{\upr((X-c)^2)+\varepsilon^2}.
\end{eqnarray}
\item[(b)]
Let $\lpr$ be defined on $\dset_{\geq}(X)=\{(X\geq c+\varepsilon), -(X-c)^2\}$, and let $\mset_{\geq}(X)=\{P: P \text{ dF-coherent}, P\geq\lpr \text{ on } \dset_{\geq}(X)\text{ and }P(X)=c\}$.
If $\mset_{\geq}(X)\neq\varnothing$, then
\begin{eqnarray}
	\label{eq:cantelli_general_2}
	\lpr(X\geq c+\varepsilon)\leq\frac{\upr((X-c)^2)}{\upr((X-c)^2)+\varepsilon^2}.
\end{eqnarray}
\end{itemize}
\end{corollary}
\begin{proof}
Apply Proposition \ref{pro:cantelli_imprecise} to $Y=X-c$ (noting, for (a), that $\dset_{\leq}(X)=\dset^{\leq}_{0}(Y)=\{(Y\leq-\varepsilon), -Y^2\}$ and $\mset_{\leq}(X)=\mset_{0}^{\leq}(Y)$; similarly for (b)).
\end{proof}
A natural choice for $c$ is $c=\lpr(X)$.
With this, \eqref{eq:cantelli_general_1} reads
\begin{eqnarray}
\label{eq:cantelli_general_1_lpr}
\lpr(X\leq\lpr(X)-\varepsilon)\leq\frac{\upr((X-\lpr(X))^2)}{\upr((X-\lpr(X))^2)+\varepsilon^2},
\end{eqnarray}
giving a non-trivial upper bound to the lower probability that $X$ is `too far away' (from below) from $\lpr(X)$.
For instance, with $\varepsilon=3\sqrt{\upr((X-\lpr(X))^2)}$ we obtain the bound
\begin{equation*}
\lpr\left(X\leq\lpr(X)-3\sqrt{\upr((X-\lpr(X))^2)}\right)\leq\frac{1}{10}.
\end{equation*}
Appealing alternative choices for $c$ are $c=\upr(X)$, $c=\frac{\lpr(X)+\upr(X)}{2}$. Note that (a) in Corollary \ref{cor:cantelli_general} does not apply when $c\notin [\inf X,\sup X]$, since than $\mset_{\leq}(X)=\varnothing$ ($P(X)=c$ being not dF-coherent).

Analogous considerations apply to inequality \eqref{eq:cantelli_general_2}.
Thus, Cantelli's inequalities may be viewed as concentration inequalities, holding for lower previsions not necessarily coherent.

\subsection{Cantelli's inequalities and coherence}
\label{subsec:cantelli_coherence}
The results of the previous section can be strengthened assuming coherence.

As a first step, take $c=\lpr(X)$ and suppose that $\lpr$ is coherent on $\{(X\leq\lpr(X)-\varepsilon), -(X-\lpr(X))^2, X\}$.
Then, $\mset_{\leq}(X)$ in Corollary~\ref{cor:cantelli_general} is non-empty, by the Lower Envelope Theorem (Theorem \ref{thm:lower_envelope}),
and inequality \eqref{eq:cantelli_general_1_lpr} obtains.
Quite analogously, coherence of $\lpr$ on $\{(X\geq\lpr(X)+\varepsilon), -(X-\lpr(X))^2, X\}$ is sufficient to establish the corresponding inequality
\begin{eqnarray}
	\label{eq:cantelli_general_2_lpr}
	\lpr(X\geq\lpr(X)+\varepsilon)\leq\frac{\upr((X-\lpr(X))^2)}{\upr((X-\lpr(X))^2)+\varepsilon^2}.
\end{eqnarray}
If in particular $\lpr=\upr=P$, with $P$ dF-coherent prevision, \eqref{eq:cantelli_general_1_lpr} and \eqref{eq:cantelli_general_2_lpr} specialise into
\begin{align}
\label{eq:cantelli_prec_var_1}
P(X\leq P(X)-\varepsilon)&\leq\frac{V_{P}(X)}{V_{P}(X)+\varepsilon^2}\\
\label{eq:cantelli_prec_var_2}
P(X\geq P(X)+\varepsilon)&\leq\frac{V_{P}(X)}{V_{P}(X)+\varepsilon^2}
\end{align}
where
\begin{equation*}
V_{P}(X)=P((X-P(X))^2)
\end{equation*}
is the variance of $X$ whenever $P$ is an expectation, in which case \eqref{eq:cantelli_prec_var_1}, \eqref{eq:cantelli_prec_var_2} boil down to \eqref{eq:cantelli_precise_prev_non_zero}.

However, stricter Cantelli's bounds than \eqref{eq:cantelli_general_1_lpr} and \eqref{eq:cantelli_general_2_lpr} apply with coherent lower previsions.
Prior to seeing this, let us recall from \cite[Appendix G]{walley_statistical_1991} the definitions of lower and upper variance.
\begin{definition}
\label{def:imprecise_variance}
For a given gamble $X$ and a coherent lower prevision $\lpr$, let
\begin{align}
\label{eq:lower_variance}
\lvx&=\min_{c\in\rset}\{\lpr((X-c)^2)\} \text{ (lower variance)}\\
\label{eq:upper_variance}
\uvx&=\min_{c\in\rset}\{\upr((X-c)^2)\} \text{ (upper variance)}
\end{align}
\end{definition}
It is also proven in \cite[Appendix G]{walley_statistical_1991} that the minima in \eqref{eq:lower_variance}, \eqref{eq:upper_variance} are attained by the previsions in the credal set $\mset$ of $\lpr$, $P_1$ and $P_2$ respectively, such that\footnote{
In this result, $\lpr$ is defined on the linear space $\lset$ of all gambles.
In the sequel we shall operate with restrictions of $\lpr$ onto the subset(s) of $\lset$ of interest.
To see that these two approaches are equivalent, recall from Section~\ref{subsec:imprecise_previsions} that given a coherent $\lpr$ on any such subset $\dset$, $\lpr$ always has a coherent extension on any superset of $\dset$, in particular on $\lset$.
}
\begin{align}
\label{eq:lower_variance_min}
\lvx&=V_{P_{1}}(X)=\min_{P\in\mset}\{V_P\}\\
\label{eq:upper_variance_max}
\uvx&=V_{P_{2}}(X)=\max_{P\in\mset}\{V_P\}.
\end{align}
Among the properties of $\lvx$, $\uvx$, we recall that \cite[Appendix G]{walley_statistical_1991}
\begin{eqnarray}
\label{eq:impr_var_property_1}
0\leq\lvx\leq\uvx,\ \lvx = \uvx \Rightarrow \lpr(X)=\upr(X).
\end{eqnarray}
Given this, we have that
\begin{proposition}
\label{pro:cantelli_coherence}
Let $X$ be a gamble and $\varepsilon>0$.
\begin{itemize}
\item[(a)]
Let $\lpr$ be coherent on $\{X, -X, (X\leq\lpr(X)-\varepsilon), (X\leq\upr(X)-\varepsilon)\}$. Then,
\begin{align}
\label{eq:cantelli_coherent_leq_1}
\lpr(X\leq\upr(X)-\varepsilon)\leq\frac{\uvx}{\uvx+\varepsilon^2}\\
\label{eq:cantelli_coherent_leq_2}
\lpr(X\leq\lpr(X)-\varepsilon)\leq\frac{\lvx}{\lvx+\varepsilon^2}.
\end{align}
\item[(b)]
Let $\lpr$ be coherent on $\{X, -X, (X\geq\lpr(X)+\varepsilon), (X\geq\upr(X)+\varepsilon)\}$. Then,
\begin{align}
	\label{eq:cantelli_coherent_geq_1}
	\lpr(X\geq\lpr(X)+\varepsilon)\leq\frac{\uvx}{\uvx+\varepsilon^2}\\
	\label{eq:cantelli_coherent_geq_2}
	\lpr(X\geq\upr(X)+\varepsilon)\leq\frac{\lvx}{\lvx+\varepsilon^2}.
\end{align}
\end{itemize}
\end{proposition}
\begin{proof}
We prove only part (a), the proof of (b) being analogous.
We prove first \eqref{eq:cantelli_coherent_leq_1}.
Because $\lpr$ is coherent, there is $P^*\in\mset$ such that $P^*\geq\lpr$ and $P^*(-X)=\lpr(-X)=-\upr(X)$, implying $P^*(X)=\upr(X)$.
Hence,
\begin{eqnarray*}
\lpr(X\leq\upr(X)-\varepsilon)\leq P^*(X\leq P^*(X)-\varepsilon)\leq\frac{V_{P^*}(X)}{V_{P^*}(X)+\varepsilon^2}\leq\frac{\uvx}{\uvx+\varepsilon^2},
\end{eqnarray*}
where the second inequality is due to \eqref{eq:cantelli_prec_var_1}, the third to \eqref{eq:upper_variance_max}.

To prove \eqref{eq:cantelli_coherent_leq_2}, let $P_1$ be the prevision in $\mset$ satisfying \eqref{eq:lower_variance_min}.
Then we have, by \eqref{eq:cantelli_prec_var_1} at the second inequality and by the properties of $P_1$ elsewhere
\begin{eqnarray*}
\lpr(X\leq\lpr(X)-\varepsilon)\leq P_1(X\leq P_{1}(X)-\varepsilon)\leq\frac{V_{P_{1}}(X)}{V_{P_{1}}(X)+\varepsilon^2}=\frac{\lvx}{\lvx+\varepsilon^2}.
\end{eqnarray*}
\end{proof}
Let us discuss the role of the new Cantelli's inequalities in Proposition~\ref{pro:cantelli_coherence}.
Considering first \eqref{eq:cantelli_coherent_leq_1} and \eqref{eq:cantelli_coherent_leq_2}, both give a bound to the lower probability $\lpr(X\leq c)$.
To apply \eqref{eq:cantelli_coherent_leq_1}, it is necessary that $c<\upr(X)$.
On the other hand, if $c<\lpr(X)$, inequality \eqref{eq:cantelli_coherent_leq_2} obtains too and is \emph{definitely preferable}.
In fact, recalling \eqref{eq:impr_var_property_1}, we have that
\begin{eqnarray}
\label{eq:rel_dx_low_upp}
\frac{\lvx}{\lvx+\varepsilon^2}\leq\frac{\uvx}{\uvx+\varepsilon^2},
\end{eqnarray}
so that the upper bound \eqref{eq:cantelli_coherent_leq_1} is looser than \eqref{eq:cantelli_coherent_leq_2}.

Inequality \eqref{eq:cantelli_coherent_leq_2} is  also weakly stricter than inequality \eqref{eq:cantelli_general_1_lpr}:
this is easily seen from \eqref{eq:rel_dx_low_upp} and \eqref{eq:upper_variance}.
However, \eqref{eq:cantelli_coherent_leq_2} is \emph{typically stricter} than \eqref{eq:cantelli_general_1_lpr},
meaning that $\lpr(X)=\upr(X)$, i.e. $X$ is given a precise evaluation, when both inequalities return the same bound.
In fact, in this case $\lvx=\upr((X-\lpr(X))^2)$, implying $\lvx=\uvx$ (recalling that $\lvx\leq\uvx\leq\upr((X-\lpr(X))^2
)$ by \eqref{eq:upper_variance} and \eqref{eq:impr_var_property_1}). Therefore, $\lpr(X)=\upr(X)$ by \eqref{eq:impr_var_property_1}.
Thus, \eqref{eq:cantelli_coherent_leq_2} is the strictest left-sided Cantelli's inequality we derived.
In a comparison with inequality \eqref{eq:cantelli_general_1_lpr}, it has to be said that \eqref{eq:cantelli_coherent_leq_2} requires the additional task of determining $\lvx$.
Ways to compute $\lvx$ and $\uvx$ are discussed in \cite[Appendix G]{walley_statistical_1991};
it is reported there that the computation of $\lvx$ is simpler than that of $\uvx$.

Turning to inequality \eqref{eq:cantelli_coherent_geq_2}, it is the strictest Cantelli's bound to $\lpr(X\geq c)$.
It works when $c>\upr(X)$ and is preferable to the looser inequalities \eqref{eq:cantelli_coherent_geq_1} and \eqref{eq:cantelli_general_2_lpr}.
On the other hand, \eqref{eq:cantelli_coherent_geq_1} and \eqref{eq:cantelli_general_2_lpr} are applicable when $c\in ]\lpr(X),\upr(X)[$ and an advantage of \eqref{eq:cantelli_general_2_lpr} is that it does not require knowing $\lvx$, nor $\uvx$ (and not even coherence).

Cantelli's inequalities assume an only limited knowledge of $X$.
As such, they cannot be expected to produce very strict bounds, in general.
To get some insight on this, a partial comparison between Cantelli's and Markov's inequalities is possible.
Precisely, let us compare the Lower Markov's inequality \eqref{eq:lower_Markov_inequality} with Cantelli's inequality \eqref{eq:cantelli_coherent_geq_2}.
Both inequalities apply if $\lpr$ is coherent, $X\geq 0$, $\varepsilon>0$ and $a=\upr(X)+\varepsilon$, so that Markov's inequality is
\begin{eqnarray}
\label{eq:Markov_inequality_comparison}
\lpr(X\geq\upr(X)+\varepsilon)\leq\frac{\lpr(X)}{\upr(X)+\varepsilon}.
\end{eqnarray}
By \eqref{eq:cantelli_coherent_geq_2}, $\lpr(X\geq\upr(X)+\varepsilon)$ is majorised by $\frac{\lvx}{\lvx+\varepsilon^2}$,
thus Cantelli's bound is better than Markov's iff
\begin{eqnarray}
\label{eq:Markov_Cantelli_comparison}
\frac{\lvx}{\lvx+\varepsilon^2}\leq\frac{\lpr(X)}{\upr(X)+\varepsilon}.
\end{eqnarray}
Discarding the extreme case $\lpr(X)=0$ (where Markov's inequality \eqref{eq:Markov_inequality_comparison} cannot be improved, implying $\lpr(X\geq\upr(X)+\varepsilon)=0$),
condition \eqref{eq:Markov_Cantelli_comparison} is equivalent to
\begin{equation*}
\lpr(X)\varepsilon^2-\lvx\varepsilon+\lvx[\lpr(X)-\upr(X)]\geq 0,
\end{equation*}
a second order inequality with
\begin{equation*}
\Delta=\lvx[\lvx+4\lpr(X)(\upr(X)-\lpr(X))]\geq[\lvx]^2\geq 0.
\end{equation*}
It ensues that
\begin{eqnarray*}
\varepsilon_1&=\frac{\lvx-\sqrt{\Delta}}{2\lpr(X)}\leq 0, \text{ while}\\
\varepsilon_2&=\frac{\lvx+\sqrt{\Delta}}{2\lpr(X)}\geq\frac{\lvx}{\lpr(X)}\geq 0.
\end{eqnarray*}
Since $\varepsilon>0$, we may conclude that Cantelli's inequality is stricter than Markov's iff $\varepsilon>\varepsilon_2$.

Note also that a simple sufficient condition for Markov's inequality to be preferable is that $\varepsilon\lpr(X)<\lvx$.

We only mention that Markov's inequality \eqref{eq:lower_Markov_inequality}, taking $a=\lpr(X)+\varepsilon$, can be confronted also with inequality \eqref{eq:cantelli_general_2_lpr}.
The procedure is similar, and again the result depends on the value of $\varepsilon$.

Finally, we point out that further Cantelli-like inequalities are easily derived from the previous ones via conjugacy.
To exemplify, the bound in \eqref{eq:cantelli_coherent_leq_2} applies by monotonicity also to $\lpr(X<\lpr(X)-\varepsilon)$, which implies that
\begin{equation*}
\upr(X\geq\lpr(X)-\varepsilon)=1-\lpr(X<\lpr(X)-\varepsilon)\geq 1-\frac{\lvx}{\lvx+\varepsilon^2}=\frac{\varepsilon^2}{\lvx+\varepsilon^2}.
\end{equation*}

\section{Conclusions}
\label{sec:conclusions}
In this paper we explored how some basic probability inequalities are modified when our uncertainty evaluations are imprecise.
It turns out that the inequalities have a certain \emph{robustness} with respect to the quality of our information.
While in classical probability theory this means that a very limited (if any) distributional knowledge of a gamble $X$ is required to apply the inequalities,
here in addition we may have different degrees of consistency for our evaluations.
As we have seen, coherence is anyway preferable.
It may ensure inferences in a larger number of situations (Section~\ref{subsec:Basic inferences Jensen}, Example~\ref{ex:Markov}) or better bounds (Section~\ref{subsec:cantelli_coherence}).
It remains to be explored whether some of these bounds can be further improved, in case referring to specific coherent models.
It is on the other hand remarkable that versions of the investigated inequalities may be employed with much weaker consistency conditions, such as $2$-coherence.

We expect that similar conclusions may be drawn with generalisations of other inequalities not explored here.
This seems to be a promising area for future work.
A question one may run into quite soon is how should concentration (or dispersion, equivalently) of lower or upper previsions be measured.
While being an issue to be deepened, we guess that there might be no unique answer.
Take for instance Cantelli's bounds.
The variance $\sigma^2(X)$ appears in the probability bounds \eqref{eq:cantelli_precise_prev_non_zero}.
It is replaced by Walley's lower variance $\lvx$ in the bound \eqref{eq:cantelli_coherent_leq_2} (or by $\uvx$ in \eqref{eq:cantelli_coherent_leq_1}), but this requires coherence.
Otherwise, the role of $\sigma^2(X)$ is taken by $\upr((X-\lpr(X))^2)$ in inequality \eqref{eq:cantelli_general_1_lpr}.
Note also that Chebyshev-like inequalities are easily obtained from Markov's inequalities \eqref{eq:lower_Markov_inequality}, \eqref{eq:upper_Markov_inequality}.
These involve different variance-like quantities.
For instance, we have that, given $b>0$,
\begin{align*}
	\upr(|X-\lpr(X)|\geq b)&\leq\frac{\upr((X-\lpr(X))^2)}{b^2}, \text{ but also}\\
	\upr(|X-\upr(X)|\geq b)&\leq\frac{\upr((X-\upr(X))^2)}{b^2}.
\end{align*}
In general, $\lvx$, $\uvx$ are coherence-dependent by construction, since they refer to the credal set $\mset$ of $\lpr$.
By the Lower Envelope Theorem (Theorem~\ref{thm:lower_envelope}), the connection between a coherent $\lpr$ and its credal set is strict, but becomes looser with previsions that avoid sure loss.
Would $\lvx$ still be a reliable variance measure in such a case?
Further, with $2$-coherence $\mset$ may even be empty.
Thus, what precisely extends variance with imprecise judgements is another question waiting for more exhaustive answers.

%

\bibliographystyle{plain}
\bibliography{PVinequalities}{}

\end{document}